\documentclass[12pt,reqno]{amsproc} %

\usepackage{amsthm}
\usepackage{amsmath}
\usepackage{amssymb}

{ \theoremstyle{plain}
\newtheorem{theorem}{Theorem}%[section]
\newtheorem{lemma}{Lemma}%[section]
\newtheorem{proposition}{Proposition}%[section]
}

{ \theoremstyle{definition}
}

\theoremstyle{plain}
\theoremstyle{remark}
\newtheorem{remark}{Remark}%[section]

\theoremstyle{remark}
\newtheorem{example}{Example}%[section]

\def\Com#1{\mathbb{C}^{#1}}
\def\R#1{\mathbb{R}^{#1}}

\def\V{\mathcal{V}}
\def\DIS{\mathcal{D}}
\def\Pos#1{{ \mathcal{K}}_{#1}}

\def\x{\overline{x}}
\def\ov#1{\overline{#1}}

 \DeclareMathOperator{\Res}{Res}
\DeclareMathOperator{\interior}{int}
\DeclareMathOperator{\rank}{rank} \DeclareMathOperator{\Dis}{Dis}
\DeclareMathOperator{\re}{Re} \DeclareMathOperator{\im}{Im}
\begin{document}

\addtolength{\baselineskip}{0.5mm}
%\begin{large}

\title[Positive trigonometric polynomials]{Positive trigonometric polynomials}

\author{Tkachev V.G.}

\address{Volgograd State University}

\email{Vladimir.Tkachev@volsu.ru}

\thanks{Paper is supported by Russian President grant for
Young Doctorates, Ref. number 00-15-99274 }



%\date{}

%%%%%%%%%%%%%%%%%%%%%%%%%%%%%%%%%%%%5%
\vspace*{-1cm}

\begin{abstract}
We study the boundary of the nonnegative trigonometric polynomials
from the algebraic point of view. In particularly, we show that it
is a subset of an irreducible algebraic hypersurface and
established its explicit form in terms of resultants.
\end{abstract}

\maketitle
%%%%%%%%%%%%%%%%%%%%%%%%%%%%%

\section*{Introduction}

Let $n\geq 1$ be an integer and
\begin{equation}
T(t)=\re \sum_{k=0}^{n}y_ke^{ikt}=\sum_{k=0}^{n}(a_k\cos
kt+b_k\sin kt) \label{eq-1}
\end{equation}
be a trigonometric polynomial of degree $n$ with real coefficients
$a_k=\re y_k$, $b_k=-\im y_k$ where $b_0=0$ and $y_0= a_0$ is a
real positive number. We say that $T(t)$ is \textit{nonnegative}
if the inequality $T(t)\geq 0$ holds for all $t\in \R{}$.

The problem of deciding whether given a trigonometric polynomial
is nonnegative has long and reach history. We refer to recent
papers \cite{GH_month}, \cite{DM} for more details and
circumstance. We mention that various reasons for the interest in
the problem of constructing nonnegative trigonometric polynomials
are: the Gibbs phenomenon, univalent functions and polynomials,
positive Jacobi polynomial sums \cite{Ask}, orthogonal polynomials
on the unit circle \cite{WAss}, zero-free regions for the Riemann
zeta-function \cite{Arest92}, \cite{ArestKon} just to mention a
few.

The well-known fact of complex analysis \cite{Akh} states that an
equivalent criterion is to check that the Hermitian Toeplitz
matrices $T_l=[c_{i-j};\; 0\leq i,j\leq l;\; c_l=0 \; \mathrm{for}
\; l>k]$ are nonnegative definite for all positive integers $l$
(here $c_k=y_k$ and $c_{-k}=\ov{y}_k$ for all $k\geq 0$). But this
condition involves infinitely many inequalities and does not give
exact information about the structure of nonnegative polynomials.

Our interest in this subject comes from the theory of univalent
algebraic polynomials in the unit disk. The simplest class of the
univalent polynomials is the body of all starlike polynomials in
the unit disk $U$. It is well-known fact that an algebraic
polynomial $P(z)$ with $P(0)=0$ is star-like univalent if and only
if $\re{} (P'(z)z/P(z))\geq 0$ in $U$. The last can be
reformulated as $P(z)/z$ has no roots in $U$ and the trigonometric
polynomial $\re{} T_P(e^{it})$ where
$$
T_P(z)=P'(z)\frac{\ov{P}(z)}{z}
$$
is nonnegative. So, an arbitrary starlike polynomial one can
associate with a certain positive trigonometric polynomial.

The main purpose of this article is translation the initial
problem into the algebraic context and our point of view of the
problem is studying of the image of $\R{+}\times \Com{n}$ under
the associated polymorphism $P\to T_P$. The associated morphism is
the special quadratic and maps $\R{+}\times \Com{n}$ onto certain
cone in the same target space.

This problem was initially studied in the lower dimensional case
by Brannan \cite{Br67}--\cite{BrBrik}. In particularly, it was
proved that the boundary of three dimensional body (the case $\deg
P=3$) is an real algebraic manifold of degree $d(3)\leq 28$.
Brannan also established the explicit representation of the
boundary of the body.

In his paper \cite{Quine} Quine was showed that the boundary of
the body of all univalent polynomials is also an real algebraic
manifold of degree at most $q(n)=4(2n^2-4n+1)(n-1)$. We notice
that this estimate is too far from to be sharp. Even in the
Brannan's case $n=3$ we have the value $q(3)=56>d(3)$.

%Further we study the boundary of nonnegative trigonometric
%polynomials cone in terms certain algebraic form.

\textsl{Acknowledgement}: We wish to thank all of staff of
Mittag-Leffler Institute for heart hospitality during the visit of
the author Mittag-Leffler Institute in 2001.

\section{Quadratic map}

\subsection{}
Let $T(t)$ is given by the representation (\ref{eq-1}). Then by
the well-known result due to L.~Fej\'er \cite{Fejer} $T(t)$ is a
nonnegative trigonometric polynomial if and only if there exists
an algebraic polynomial $X(z)=x_0+x_1z+\ldots +x_nz^n$ with
complex coefficients $x_k$ such that
\begin{equation}\label{eq-2}
T(t)=\left|X(e^{it})\right|^2.
\end{equation}
Moreover, such representation with normalization $x_0>0$ is
unique. In sequel we fix this normalization.

But still this characterization does not make it easy to decide if
a given trigonometric polynomial is nonnegative. We consider the
representation in more details.

We use the following notations. Given an algebraic polynomial
$P(z)$ of the $n$-th degree let us denote by $P^{*}(z)$ its
reciprocal double
$$
P^{*}(z)\equiv z^{n}\overline{P}(1/z),
$$
where $\overline{P}(\zeta)=\overline{P(\overline{\zeta})}$ is the
conjugate polynomial to $P(z)$. By $\Res{ (P;Q)}$ we denote the
resultant of two corresponding polynomials. We also usually
identify a polynomial $P$ and the vector in $\R{+}\times\Com{n}$
(or in $\R{2n+1}$) which consists of the coefficients of $P$.

Now we can reformulate the initial problem (\ref{eq-2}) in more
appropriate terms. Namely, we associate the trigonometric
polynomial $T(t)$ in (\ref{eq-1}) with the unique algebraic
polynomial $Y(z)$ of the form $\sum_{k=0}^{k=n}{y_kz^k}$ with
$y_0$ to be a positive real. Then the problem (\ref{eq-2}) is
equivalent to following: \textsl{Given a vector $y\in \R{+}\times
\Com{n}$ find a vector $x\in \R{+}\times \Com{n}$ such that}
\begin{equation}\label{eq-3}
\begin{split}
&y_0=2\Phi_0(x_0,x_1,\x_1,\ldots x_n,\x_n)\equiv\frac{1}{2}
\sum_{j=0}^{n} \ov{x}_jx_j\\
&y_m=\Phi_m(x_0,x_1,\x_1,\ldots x_n,\x_n)\equiv
\sum_{k=0}^{n-m}\overline{x}_kx_{k+m}, \qquad m=0,1,\ldots, n.
\end{split}
\end{equation}
If such a vector $x$ exists then we say that $y$ is a
\textit{nonnegative} vector of $\R{+}\times \Com{n}$.

It follows from formulae (\ref{eq-3}) that
\begin{equation}\label{eq-5}
X(z)X^{*}(z)=\ov{\Phi}_n(x)+\ov{\Phi}_{n-1}(x)z+\ldots+2\Phi_{0}(x)z^{n}+
\Phi_{1}(x)z^{n+1}+\ldots+\Phi_{n}(x)z^{2n},
\end{equation}
or in the other form,
\begin{equation}\label{eq-31}
X(z)X^{*}(z)=Y^\star(z)+z^nY(z).
\end{equation}

\begin{remark}
\label{rem-1} The mentioned above problem has a clear geometric
interpretation. Really, the set of all nonnegative vectors forms a
cone $\Pos{n}$ in $\R{+}\times \Com{n}\equiv \R{2n+1}$ which is
strongly contained in $\R{2n+1}$. Moreover, it follows that a
point $y$ is an inner point of $\Pos{n}$ if and only if there
exists a point $x\in \R{+}\times \Com{n}$ such that the mapping
$\Phi$ has the highest rank at $x$. Other words,
\begin{equation}\label{rang}
y\in \interior{\Pos{n}} \Leftrightarrow \max_{x\in
\Phi^{-1}(y)}\rank d_{\R{}}\Phi(x)=2n+1.
\end{equation}
As a consequence, the boundary $\partial \Pos{n}$ is contained in
the image $\Phi(\V_n)$ where $\V_n$ is zero-set of the Jacobian
$\det d_{\R{}}\Phi(x)$. Nevertheless, in general case only the
strong inclusion $\partial \Pos{n}\varsubsetneqq\Phi(\V_n)$ can be
possible.
\end{remark}

\subsection{}
Our first step is to describe the structure of the image
$\Phi(V_n)$. To do this we need some auxiliary assertions.

\begin{lemma}
We have
\begin{equation}\label{eq-4}
\det d_{\R{}}\Phi(x)\equiv \det \frac{\partial
(\Phi_0,\Phi_1,\ov{\Phi}_1,\ldots,\Phi_n,\ov{\Phi}_n)}{\partial
(x_0,x_1,\ov{x}_1,\ldots,x_n,\ov{x}_n)}=2x_0\Res{(X,X^{*})}.
\end{equation}
\label{lem-1}
\end{lemma}

\begin{proof}
Let us denote by $R(z)$ the polynomial $X(z)X^{*}(z)$. Then we
have from (\ref{eq-3}) the expressions for the derivatives
\begin{equation}\label{eq-6}
\begin{split}
&R_j\equiv \frac{\partial R}{x_j}=\frac{\partial X}{x_j}\;
X^{{*}}=z^{j}X^{*},\\
&R_{-j}\equiv \frac{\partial R}{\ov{x}_j}=X\;\frac{\partial
X^{*}}{\ov{x}_j}=z^{n-j}X,\\
&R_0\equiv \frac{\partial R}{x_0}=\frac{\partial X}{x_0}\;
X^{{*}}+X\;\frac{\partial X^{*}}{x_0}=X^{*}+z^nX.
\end{split}
\end{equation}

On the other hand, the coefficients of polynomials $R_j$, $R_{-j}$
and $R_0$ produce the corresponding strings of the Jacobian
$d_{\R{}}\Phi(x)$. Moreover, we claim  that
\begin{equation}\label{eq-7}
X(z)(z^n-Q_{n-1}(z))=P_n(z)X^{*}(z),
\end{equation}
where $Q_{n-1}$ and $P_n$ are some polynomials of degrees $(n-1)$
and $n$ respectively. Indeed, we can assume that
$$
Q_{n-1}\equiv z^n-\frac{1}{x_0}X^{*}(z),\qquad P_n\equiv
\frac{1}{x_0}X(z).
$$

Moreover, the last expressions imply that $P_n(0)=X(0)/x_0=1$ and
it follows from (\ref{eq-7}) that
\begin{equation}\label{eq-8}
X(z)z^n=\biggl(q_0+q_1z+\ldots
+q_{n-1}z^{n-1}\biggr)X(z)+\biggl(1+p_1z+\ldots
+p_{n}z^{n}\biggr)X^{*}(z).
\end{equation}
Here $q_k$ and $p_k$ are the coefficients of the corresponding
capitals polynomials.

By virtue of (\ref{eq-6}) we obtain that
\begin{equation}
X(z)z^n-X^{*}(z)=q_0R_{-n}+q_1R_{-n+1}+\ldots
+q_{n-1}R_{-1}+p_1R_{1}+\ldots +p_{n}R_{n},
\label{eq-91}
\end{equation}
and
$$
2X(z)z^n\equiv
\biggl(X(z)z^n-X^{*}(z)\biggr)+\biggl(X(z)z^n+X^{*}(z)\biggr)=
\biggl(X(z)z^n-X^{*}(z)\biggr)+R_0(z).
$$

By (\ref{eq-91}) the last equality can be interpreted as the fact
that $2X(z)z^n$ is a linear combination of $R_0$ and other strings
of the Jacobian matrix $d_{\R{}}\Phi(x)$ (here and further on we
use the mentioned above equivalence between the polynomial and
vectors).

Now we recall that by Sylvester's formula for the resultant (see
\cite{WanDer})
\begin{multline}\label{eq-9}
\det d_{\R{}}\Phi(x)= 2\det[z^nX,z^{n-1}X^{*},z^{n-1}X,\ldots,
z^{j}X^{*},z^{n-j}X,\ldots, z^{n}X^{*}, X]=\\
=2(-1)^{\sigma}\det[X,z X,\ldots, z^{n}X,X^{*},zX^{*},\ldots,
z^{n}X^ {*}]=2(-1)^{\sigma}\Res{(zX^{*},X)}.
\end{multline}
Here by $\det [h_1,\ldots,h_k]$ we denote the determinant of the
matrix with strings $h_k$. The symbol $\sigma$ is equal to the
corresponding permutation number between to matrices. It is easy
to compute that $\sigma =n^2+n\equiv 0 \mod 2$ and by the
multiplication property of resultant we arrive at
$$
\det
d_{\R{}}\Phi(x)=2\Res{(zX^{*},X)}=2\Res{(z,X)}\Res{(X^{*},X)}=2x_0\Res{(X^{*},X)}
$$
and the lemma is proved completely.
\end{proof}

\subsection{}
We introduce the special notation for the characteristic
$\Res(X,X^{*})$. Namely, we call this quantity the
\textit{M\"obius discriminant} of $X(z)$ and will denote it by
$V_n(X)$. The direct consequence of its definition and Sylvester's
formula is that $V_n(X)$ is a homogeneous form of $X$ (or it the
same as the form of variables $x_k$ and $\ov{x}_k$) of degree
$\deg V_n=2n$. Moreover,
$$
V_n(X)=(-1)^nV_n(X^{*}),
$$
because $\Res(P,Q)=(-1)^{\deg P\deg Q}\Res(Q,P)$

Further we use the fact that $V_n(X)$ is actually a real valued
form. To see that we notice that by the decomposition property of
resultant and relation between the roots of $X$ and $X^{*}$
\begin{equation}\label{eq-10}
V_n(X)=|x_n|^{2n}\prod_{j,k=1}^{n}(z_j\ov{z}_k-1).
\end{equation}
where $z_j$  are the roots of $X(z)$.

Really, we recall that for arbitrary polynomials $P(z)$ and $Q(z)$
of $n$ and $m$ degrees one holds the Cayley formula for its
resultant
\begin{equation}
\Res{(P;Q)}=p_n^mq_m^n\prod_{j=1}^{n}\prod_{k=1}^{m}(a_j-b_k),
\label{eq-11} \end{equation} where $a_j$ and $b_k$ are the roots
of $P(z)$ and $Q(z)$ respectively.

So, if we have decomposition $ X(z)=x_0+x_1z+\ldots
+x_nz^n=x_n\prod_{j=1}^{n}(z-z_j),
$
then by the definition of $X^{*}$
$$
X^{*}(z)=\ov{x}_n+\ov{x}_{n-1}z+\ldots
+x_0z^n=z^n\ov{X}(1/z)=\ov{x}_n\prod_{j=1}^{n}(1-z\ov{z}_j).
$$
It follows from (\ref{eq-11}) and  Vi\`ete formula that
\begin{equation}
\label{eq-12}
V(X)\equiv\Res{(X;X^{*})}=x_n^nx_0^n\prod_{j=1}^{n}\prod_{k=1}^{n}(z_j-1/\ov{z}_k)=
|x_n|^{2n}\prod_{j,k=1}^{n}(z_j\ov{z}_k-1),
\end{equation}
and (\ref{eq-10}) is proved.

\subsection{}
To obtain the explicit representation for the boundary $\partial
\Pos{n}$ we use arguments of Remark~\ref{rem-1}. We need a
suitable generalization of discriminant. Let $Y(z)$ be an
arbitrary polynomial of degree $n$ with real initial coefficient
$y_0$. Let $\Dis(Y)$ is the discriminant of $Y$ (see
\cite[\S~33]{WanDer}) and
$$
\Dis_2(Y)=\Dis(Y^{*}(z)+z^nY(z)).
$$
The last expression can be regarded as certain $A$-discriminant in
terminology of  \cite[c.~271]{GKZ}. It is the direct corollary of
homogeneity of $\Dis (Y)$  that the 2nd discriminant $\Dis_2(Y)$
has $(4\deg Y-2)$ degree.

\begin{lemma}
\label{lem-2} If \ $Y=\Phi(X)$ then
\begin{equation}\label{eq-13}
\Dis_2(Y)=|\Dis(X)|^2\,V^2(X).
\end{equation}
\end{lemma}

\begin{proof}
First we notice that by the discriminant multiplicative rule (see.
\cite[Theorem 3.4]{Pra}) we have
\begin{equation}\label{eq-121}
\Dis(PQ)=\Dis(P)\Dis(Q)\Res^2(P;Q),
\end{equation}
which yields from the definition of the M\"obius discriminant and
identity  (\ref{eq-31}) that
$$
\Dis_2(Y)=\Dis(X)\Dis(X^\star)\,V^2(X).
$$
On the other, it is easy to see by the definition that the
reciprocal polynomial $X^\star$ has conjugate to $X$ discriminant
which completes  the proof.
\end{proof}

We denote by $\V_n$ and $\DIS_n$ the zero-sets $V(X)=0$ and
$\Dis(X)=0$ respectively. Then the following assertion shows the
shadow-character property of $\DIS_n$ with respect to mapping
$\Phi$.

\begin{lemma}\label{lem-3}
We have $\Phi(\DIS_n)\subset \Phi(\V_n)$.
\end{lemma}

\begin{proof}
Due to (\ref{eq-5}) it is sufficient to prove that given a
polynomial $X\in\DIS_n$ there exists a polynomial $Q(z)$ with
positive initial coefficient which is in $\V_n$ and the equality
$XX^{*}=QQ^{*}$ holds.

Let $X$ be arbitrary polynomial of $n$-th degree regarded as a
point in $\DIS_n$. Then $\Dis(X)=0$. It follows from the
definition of discriminant that there exist a multiple root $z_p$
of $X$. We can assume that $p=n$ and therefore
$$
X(z)=x_n (z-z_n)^2\prod_{j=1}^{n-2}(z-z_j).
$$

We consider the following polynomial
$$
Q(z)=\lambda x_n
(z-z_n)(z-1/\ov{z}_n)\prod_{j=1}^{n-2}(z-z_j)=q_0+q_1z+\ldots
+q_nz^n,
$$
where $\lambda$ is some positive  real number which is specified
later. Then from its definition $V(Q)=0$ and hence $Q\in \V_n$.
Moreover,
$$
q_0=Q(0)=\lambda X(0)/|z_n|^2=\lambda x_0/|z_n|^2
$$
is a positive real and therefore
$Q(z)$ is admissible.

On the other hand
$$
Q^{*}(z)=z^n \ov{Q}(1/z)=\lambda \ov{x}_n
(1-z\ov{z}_n)(1-z/z_n)\prod_{j=1}^{n-2}(1-z\ov{z}_j),
$$
and it is easy to compute that
$$
Q(z)Q^{*}(z)=\frac{\lambda^2}{|z_n|^2}X(z)X^{*}(z).
$$

Consequently, the choice $\lambda =|z_n|$ lead us to required
equality $QQ^{*}=XX^{*}$ and the lemma is proved.
\end{proof}

Now the main result of this section follows from (\ref{rang}) and
previous lemmas:

\begin{theorem}\label{theo-2} The set $\Pos{n}$ of nonnegative
polynomials $Y$ is a closed cone in $\R{2n+1}$ (with nonempty
interior) which boundary $\partial \Pos{n}$ is contained in the
algebraic hypersurface $\Dis_2(Y)=0$.
\end{theorem}

We emphasize that the theorem gives the sharp description of
$\Pos{n}$ because the boundary $\Pos{n}$ turns out to be a part of
the irreducible hypersurface (see the next section).

\begin{example}\label{ex-1}
Let $n=1$ and we use the previous notations. Then it is easy to
compute that $V(X)\equiv V(x_0,x_1)=x_0^2-|x_1|^2$. In this case
we have the very similar expression $\Dis_2(Y)\equiv
\Dis_2(y_0,y_1)=4(y_0^2-|y_1|^2)$.
\end{example}

\begin{example}\label{ex-2}
In the two-dimensional case $n=2$ we have
$$
V(x_0,x_1,x_2)=(x_0^2-|x_1|^2)^2-(x_0x_1-\ov{x}_1x_2)(x_0\ov{x}_1-\ov{x}_2x_1),
$$
and
\begin{equation*}
\begin{split}
 \Dis_2(Y)= &36y_1^3 \ov{y}_1 y_0 \ov{y}_2-320 y_2 y_1
\ov{y}_1 y_0^2 \ov{y}_2- 4 y_1^3 \ov{y}_1^3+256 y_2 y_0^4
\ov{y}_2-32 y_2 y_0^3 \ov{y}_1^2- 27 y_1^4 \ov{y}_2^2- \\
        & -512 y_2^2\ov{y}_2^2 y_0^2+288 y_2^2 \ov{y}_2 y_0 \ov{y}_1^2+
        36 y_2y_1\ov{y}_1^3y_0+4 y_0^2 y_1^2 \ov{y}_1^2- 32 y_0^3 y_1^2 \ov{y}_2-\\
        &-192 y_2^2
\ov{y}_2^2 y_1 \ov{y}_1- 6 y_2 \ov{y}_2 y_1^2 \ov{y}_1^2+288 y_2
\ov{y}_2^2 y_0 y_1^2+ 256 y_2^3 \ov{y}_2^3-27 y_2^2 \ov{y}_1^4
\end{split}
\end{equation*}
\end{example}

\section{Irreducibility}

\subsection{}
The aim of this paragraph is to establish irreducibility of
$\Dis_2(Y)$ form which produces the nonnegative polynomials
boundary. Let $y_j=a_j-ib_j$, $j=1,\ldots,n$ and $b_0=0$.

\begin{theorem} The form
$\Dis_2(Y)=\Dis_2(y_0,y_1,\ov{y}_1,\ldots,y_n,\ov{y}_n)$ as well
as its real representee
$\Dis_{2,\R{}}(a_0,a_1,b_1,\ldots,a_n,b_n)$ are irreducible over
$\Com{}$. \label{theo-3}\end{theorem}

The following auxiliary assertion seems to be known but we can not
arrange any citation on it in the wide literature and because we
give its proof for completeness.

\begin{lemma}
The form $\Dis{}(p_0,p_1,\ldots,p_m)\equiv\Dis{}(P)$ is
irreducible in $\Com{}[p_0,p_1,\ldots,p_m]$. \label{lem-4}
\end{lemma}
\begin{proof}
Let us suppose that there exists nontrivial factorization
\begin{equation}
\Dis{}(P)=A(P)B(P) \label{eq-14}
\end{equation}
where $A(P)$ and $B(P)$ are different from constants. Because
$\Dis{} (P)$ is homogeneous then $A(P)$ and $B(P)$ are homogeneous
too. Moreover, $\min \{\deg_p A,\deg_p B\}\geq 1$, where $\deg_p$
is the degree with respect to $\Com{}[p_0,p_1,\ldots,p_m]$.

We consider factorization $P(z)=p_m\prod _{k=1}^m(z-z_i)$ where
$z_i$ are the roots of $P(z)$ with their multiplicity. Then every
fraction $p_k/p_m$, $k=1,\ldots, m-1$ is the $k$-th elementary
symmetric function of $z_i$. Other words,
\begin{equation}\label{eq-15}
    (-1)^k\frac{p_k}{p_m}=\sum_{1\leq i_1<\ldots <i_k\leq
    m}z_{i_1}\ldots z_{i_k}\equiv \sigma _k(z_1,\ldots,z_m),
    \end{equation}
and it follows from  (\ref{eq-15}) that
\begin{equation}\label{eq-16}
\begin{split}
&A(P)=p_m^{\deg_p A}A(\sigma_1, \ldots ,\sigma_{m-1},1)\equiv
    p_m^{\deg_p A}A_1(z_1, \ldots ,z_m), \\
& B(P)=p_m^{\deg_p B}B(\sigma_1, \ldots ,\sigma_{m-1},1)\equiv
    p_m^{\deg_p B}B_1(z_1, \ldots ,z_m).
\end{split}
\end{equation}
Obviously, both of $A_1(z_1,\ldots,z_m)$ and $B_1(z_1,\ldots,z_m)$
are symmetric polynomials by their appearance. On the other hand,
$\deg_p \Dis {}=\deg_p{}A+\deg_p{}B$ and it follows
\begin{equation}\label{eq-17}
    \Dis{} (P)=p^{2n-2}_m\prod_{m\geq i>j\geq 1}(z_i-z_j)^2.
    \end{equation}
Hence we have from (\ref{eq-16}) and (\ref{eq-17}) that
    \begin{equation}\label{eq-18}
\prod_{m\geq i>j\geq 1}(z_i-z_j)^2=
A_1(z_1,\ldots,z_m)B_1(z_1,\ldots,z_m).
\end{equation}
But the left-hand side of (\ref{eq-18}) does not admit any
nontrivial decomposition in $\Com{} [z_1,\ldots ,z_m]$ on two
symmetric polynomial multipliers (see \cite[Chapt.~V]{Litl}).

It follows that  $A_1(z_1,\ldots,z_m)$ (or $B_1(z_1,\ldots,z_m)$)
is a constant. This can occur if and only if
$A(P)=c_1p_m^{\deg_p{} P}$ for some $c_1 \ne 0$. But this means
that $p_m$ divides $\Dis {}(P)$, and, in particular,
$$
\Dis {}(p_0,p_1,\ldots , p_{m-1},0)\equiv 0.
$$
But the last property is not the case because by the elementary
discriminant property
$$
\Dis{}(p_0,p_1,\ldots,p_{m-1},0)=\pm p^2_{m-1}\Dis{}
(p_0,p_1,\ldots,p_{m-1})\ne 0.
$$
The contradiction obtained proves the lemma.
\end{proof}

%\subsection{}

\begin{proof}[Proof of Theorem~\ref{theo-3}]
First we notice that the forms are irreducible or not
simultaneously because there exists $\Com{}$-linear isomorphism
between $\{y_0,\ldots,\ov{y}_n\}$ and $\{a_0,a_1,b_1,\ldots,
b_n\}$. So it sufficient to check irreducibility of the first form
only.

We notice that $\Dis_2{}(Y)=\Dis{}(\ov{y}_n, \ov{y}_{n-1},\ldots
,\ov{y}_1,2y_0,y_1,\ldots ,y_n)$ and it follows from previous
lemma that $\Dis{} (Y)$ is irreducible in $\Com{}[y_0,y_1,\ldots,
y_n,\ov{y}_1, \ldots, \ov{y}_n]$, as well as $\Dis_2{}(Y)$ and the
theorem is proved.
\end{proof}

%\end{large}

\end{document}